\documentclass[a4paper]{article}
\usepackage{graphicx}
\usepackage{booktabs}
\usepackage{xcolor}
\usepackage{amssymb}
\usepackage{amsfonts}

\newcommand{\R}{\mathbb{R}}

\usepackage{lineno}
\usepackage{newtxmath}
\usepackage{natbib}
\usepackage{geometry}
\usepackage{comment}
\usepackage{subfiles} % Best loaded last in the preamble
\usepackage{academicons}
\definecolor{orcidlogocol}{HTML}{A6CE39}
\newcommand{\lapprox} {\, \lower3pt\hbox{$\sim$}\llap{\raise2pt\hbox{$<$}}\,}
\newcommand{\gapprox} {\, \lower3pt\hbox{$\sim$}\llap{\raise2pt\hbox{$>$}}\,}

\chardef\us=`\_

\begin{document}
\Large
\textbf{Three-dimensional numerical schemes for the segmentation of the psoas muscle in X-ray computed tomography images}

\normalsize

Giulio Paolucci$^1$, Isabella Cama$^1$, Cristina Campi$^{1,2}$, Michele Piana$^{1,2}$ \\

\hspace{-0.5cm}$^1$ MIDA, Dipartimento di Matematica, Università di Genova, via Dodecaneso 35 16146 Genova, Italy \\
email: volpara@dima.unige.it,piana@.dima.unige.it, massone@.dima.unige.it \\
$^2$ IRCCS Policlinico San Martino Genova, Largo Rossana Benzi 10 16132 Genova, Italy \\

\date{\today}

\begin{center}
   \textbf{Abstract} 
\end{center}

The analysis of the psoas muscle in morphological and functional imaging has proved to be an accurate approach to assess sarcopenia, i.e. a systemic loss of skeletal muscle mass and function that may be correlated to multifactorial etiological aspects. The inclusion of sarcopenia assessment into a radiological workflow would need the implementation of computational pipelines for image processing that guarantee segmentation reliability and a significant degree of automation. The present study utilizes three-dimensional numerical schemes for psoas segmentation in low-dose X-ray computed tomography images. Specifically, here we focused on the level set methodology and compared the performances of two standard approaches, a classical evolution model and a three-dimension geodesic model, with the performances of an original first-order modification of this latter one. The results of this analysis show that these gradient-based schemes guarantee reliability with respect to manual segmentation and that the first-order scheme requires a computational burden that is significantly smaller than the one needed by the second-order approach. 
$$ $$  
\textbf{key words.} image segmentation, X-ray Computed Tomography (CT), three-dimensional level set methods, sarcopenia

%\titlerunning{Variation of Electron Flux Spectra using STIX Observations}
%\authorrunning{Volpara et al.}

%\maketitle

\section{Introduction}\label{sec1}

Sarcopenia \citep{cruz2019sarcopenia} is a pathological condition of the skeletal muscle characterized by loss of mass, and depletion of strength and physical performance. In the last decade, several studies have shown that sarcopenia is associated to specific chronic diseases. Just as examples, and not exhaustively, low skeletal muscle mass is associated with post-transplantation death \citep{kaido2013impact,masuda2014sarcopenia}; it is a risk factor for poor performances and toxicity of chemotherapy intervention in cancer patients \citep{vergara2020sarcopenia,chindapasirt2015sarcopenia, pamoukdjian2018prevalence,collins2014assessment,villasenor2012prevalence}; it has a high probability of co-occurrence with several neurodegenerative conditions such as Parkinson's disease and multiple sclerosis \citep{drey2017associations,yuksel2022sarcopenia}.

Despite its relevance as a frailty biomarker, sarcopenia is difficult to quantitatively determine, most studies in this field currently relying on estimates of variation of the mass of the psoas muscle determined by post-processing specific slices of low resolution X-ray computed tomography (CT) scans \citep{derstine2018skeletal}. Specifically, in a typical pipeline utilizing clinical data, abdominal CT images of the patient are downloaded from the hospital Picture Archiving and Communication System (PACS); a cross-sectional image in correspondence with the third lumbar vertebra (L3) is selected; the psoas is identified in this image; and the corresponding psoas area and density are computed.

Different approaches to automated two-dimensional psoas segmentation have been proposed. For example, anatomy-based manually extracted prior shapes have been utilized to initialize energy-based iterative schemes \citep{inoue,Kamiya2012,CHEN2019120}. More recently, artificial intelligence (AI) relied on convolutional neural networks (CNNs) and Generative Adversarial Networks (GANs) to realize CT images segmentation specifically tailored to psoas identification at L3 \citep{Hashimoto2019,Duong2022}. Of course, deep learning algorithms generally designed for muscles and internal organs segmentation \citep{Kamiya2018,villarini} can be applied to the specific task of psoas identification.

This image-based approach to sarcopenia quantification relies on the controversial assumption that a single CT slice is representative of the whole psoas condition \citep{manabe2023usefulness,bauckneht2022opportunistic,zopfs2020single}. Very recently \citep{bauckneht2020spinal,bauckneht2022opportunistic}, two studies conceived in our (extended) research group accounted for this ambiguity and showed that a complete three-dimensional segmentation of the psoas muscle in CT images of oncological and neurological patients may be exploited to obtain metabolic information about the diseased tissue. Both studies utilized hybrid data made of low-dose CT and positron emission tomography (PET) data of patients injected with 18F-Fluorodeoxyglucose (FDG-PET/CT), and, in both cases, the computational paradigm was based on the following three steps \citep{sambuceti2012estimating,fiz2015allogeneic,marini2018interplay}:
\begin{enumerate}
\item A segmentation algorithm was applied to co-registered CT data in order to automatically identify the psoas muscle.
\item A binary mask was constructed, where each voxel within the image region corresponding to the muscle was given value $1$ and all other voxels were given value $0$.
\item The binary mask was voxel-wise multiplied with the FDG-PET volumes in order to extract the metabolic information associated to the psoas tissue.
\end{enumerate}
Using this scheme, paper \citep{bauckneht2022opportunistic} introduced a novel image-driven biomarker, the Attenuation Metabolic Index (AMI), that can be applied as a prognostic tool in metastatic castration-resistant prostate cancer patients, while paper \citep{bauckneht2020spinal} proved that there is a common underlying mechanism for skeletal muscle and spinal cord hypermetabolism in amyotrophic lateral sclerosis (ALS). 

Of course, the reliability of these metabolic outcomes strongly depends on the accuracy with which the segmentation step 1 above is realized. Specifically, in the first study \citep{bauckneht2022opportunistic} the segmentation approach was based on an extended version of the Hough transform \citep{beltrametti2013hough}, which can be used only if a realistic parametric model for the muscle is at disposal; the second study \citep{bauckneht2020spinal} utilized a rather heuristic approach based on histogram equalization and an $\alpha$-shape algorithm to identify the region corresponding to the inner muscle. However, despite the increasing potential relevance of this kind of analysis for diagnostic/prognostic applications, no systematic study has been performed so far on the effectiveness of three dimensional numerical schemes for image segmentation of the psoas muscle in low-resolution CT data.

The objective of the present paper is to address the computational problem of segmenting the psoas muscle in clinical low-dose CT data by means of three dimensional numerical schemes originating from the level set methodology \citep{osher2001level}. Specifically, we focused on gradient-based methods \citep{hell2015approach} and considered a classical evolution model, a standard second-order geodesic model, and a first-order modification of the geodesic model that is formulated for the first time in the present paper. In all three cases the model discretization has been performed by means of a standard finite-difference scheme. As a result, we obtained three segmentation algorithms, all with the nice property that they require the users to select just one point for each psoas in the data volume to initialize the segmentation process. The application to CT images contained in hybrid PET/CT volumes showed a notable segmentation accuracy, although with rather significant differences in the required computational burden. However, due to the low spatial resolution of these clinical CT data, in some cases the segmented shapes presented some artifacts at the end of the prescribed iterations. Therefore, we have formulated and implemented an algorithm for removing these false positives, which can be applied at the end of the evolution process when necessary. 

The plan of the paper is as follows. Section \ref{evolution_models} describes the level set models employed in the study, together with the numerical scheme utilized for their discretization. Section \ref{pre_post_processing} discusses the more heuristic approaches exploited to pre-processing the CT images and post-processing the outcomes of the numerical evolution process. Section \ref{numerical_results} compares the results provided by the different numerical procedures when applied to clinical CT data. Our conclusions are offered in Section \ref{conclusions}.

\section{Evolution models and their discretization}\label{evolution_models}
%{\bf{(GIULIO PERDONAMI MA SECONDO ME IN QUESTO PARAGRAFO CONTINUANO A ESSERCI PROBLEMI CON LE DEFINIZIONI E CON LE EQUAZIONI. INOLTRE CI SONO PROBLEMI CON LA NUMERAZIONE DI ALCUNE EQUAZIONI. PUOI PER FAVORE CONTROLLARE E CORREGGERE?)}}
The general framework considered in the present study assumes that segmentation is realized by letting a three-dimensional wave-front 
\begin{equation}\label{eq:wave-front}
\gamma_t:=\{x \in \R^3 : v(t,x)=0\}
\end{equation}
evolve from an initial point inside the shape to segment until it reaches its boundary. Following the level set method \cite{osher2001level}, the field $v=v(t,x)$ is assumed to be the solution of the Cauchy problem
\begin{equation}
\label{eq_classic}
\left\{
\begin{array}{ll}
v_t+g(x)|D v|=0,\qquad&(x,t)\in \mathbb R^3\times \mathbb (0,T)\\
v(0,x)=v_0(x)&x\in \mathbb R^3,
\end{array}
\right.
\end{equation}
where $D$ is the three-dimensional gradient operator, $g(x)$ is the speed in the normal direction and $v_0$ is a proper representation of the front $\gamma_0$ (i.e. denoting by $\Omega_0$ the region enclosed by the front, $v_0$ has to be equal to $0$ on $\gamma_0$, positive outside $\Omega_0$ and negative inside). It is well known that, under appropriate assumptions on the velocity field $g(x)$ and the initial condition $v_0$, the viscosity solution of (\ref{eq_classic}) exists and is unique (see, e.g., \cite{crandall1983viscosity}). In the following we refer to the eikonal equation \eqref{eq_classic} as the classical model.
This approach can be generalized to a second order Hamilton-Jacobi equation of the form
\begin{equation}
\label{eq_HJ}
v_t+H(x,Dv,D^2v)=0,
\end{equation}
which includes more sophisticated definitions of the velocity field $g$, possibly dependent on the curvature of the surface, as for example in the Mean Curvature Motion. Here we consider the version proposed by Caselles et al. in \cite{caselles1997geodesic} denominated geodesic model
%where, in the classical case, the Hamiltonian $H$ is the continuous function
%\begin{equation}\label{eq:hamiltonian-1}
%    H(x,Dv) =  g(x)|D v|.
%\end{equation}
%The standard second-order geodesic model assumes
\begin{equation}\label{eq:hamiltonian-2}
    H(x,Dv,D^2v) = - g(x) |D v|
\left(\varepsilon \nabla \cdot \left(\frac{Dv}{|Dv|}\right)-\mu\right) - \eta Dg \cdot Dv,
\end{equation}
where $\varepsilon, \mu, \eta$ are positive parameters to appropriately tune. We also propose a first-order modification of equation (\ref{eq:hamiltonian-2}), i.e.
\begin{equation}\label{eq:hamiltonian-3}
    H(x,Dv) = \mu g(x) |D v|- \eta Dg \cdot Dv ,
\end{equation}
which is characterized by a lower computational burden. In order to complete the description of the segmentation model we have to define the velocity $g(x)$ such that it forces the front towards the relevant boundary and then it stops its evolution. A typical definition in gradient-based methods is the following
\begin{equation}\label{vel_c1}
g(x)=\frac{1}{\left(1+|\nabla(G\ast I(x))|^p\right)},\qquad p\geq 1,
\end{equation}
where $I$ represents the intensity values of the image, and $G$ a Gaussian filter of deviation $\sigma$, used to regularize noisy images. Here we use \eqref{vel_c1} with $p=2$ as in \cite{caselles1997geodesic}. 

We now provide some details about the numerical schemes utilized for the discretization of the first order equations (\ref{eq_classic}) and (\ref{eq:hamiltonian-3}), which can also be adapted to equation (\ref{eq:hamiltonian-2}) with minor changes. For sake of simplicity the description of these schemes is provided in the two-dimensional setting, although the codes we used for data analysis have been implemented for three-dimensional segmentation.

We realized the numerical discretization of all three models by means of the standard difference form
\begin{equation}\label{eq:FD_2D}
u^{n+1}_{i,j} =  S_\Delta(u^{n})_{i,j} := u^{n}_{i,j} -\Delta t ~ h\left(x_j,y_i,D_x^-u^n_{i,j},D_x^+u^n_{i,j},D_y^-u^n_{i,j},D_y^+u^n_{i,j}\right),
\end{equation}
where the numerical Hamiltonian is the
local Lax-Friedrichs Hamiltonian
	\begin{align}
	\label{local_lax_fried_2D}
	h(x,y,p^-,p^+,q^-,q^+) := &H\left(x,y,\frac{p^++p^-}{2},\frac{q^++q^-}{2}\right)\nonumber\\
	&-\frac{\alpha_x(p^-,p^+)}{2}(p^+-p^-)-\frac{\alpha_y(q^-,q^+)}{2}(q^+-q^-),
	\end{align}
	with
	\begin{equation}
	\alpha_x(p^-,p^+) := \max_{x,y,q,\atop p\in I(p^-,p^+)}\left|H_p(x,y,p,q)\right|,\qquad \alpha_y(q^-,q^+) := \max_{x,y,p,\atop q\in I(q^-,q^+)}\left|H_q(x,y,p,q)\right|,
	\end{equation}
	with $I(a,b) := [\min(a,b),\max(a,b)]$. 
We point out that
$h(x,y,p^-,p^+,q^-,q^+)$ is a Lipschitz continuous function, with 
\begin{equation}\label{pippo-1}
  D_x^{\pm} u^n_{i,j}:=\pm \frac{u^n_{i,j\pm 1}-u^n_{i,j}}{\Delta x},
\end{equation} 
and 
\begin{equation}\label{pippo-2}
 D_y^{\pm}u^n_{i,j}:=\pm \frac{u^n_{i\pm 1,j}-u^n_{i,j}}{\Delta y}.   
\end{equation}
Further, in this context consistency is imposed via condition
\begin{equation}\label{eq:consistency}
   h(x,a,a,b,b)=H(a,a,b), 
\end{equation}
whereas monotonicity requires the numerical Hamiltonian $h$ to be non-decreasing with respect to its second and fourth argument and non-increasing with respect to the third and the fifth ones. The resulting scheme \eqref{eq:FD_2D}-\eqref{local_lax_fried_2D} is monotone under the Courant–Friedrichs–Lewy (CFL) type condition $\frac{\Delta t}{\Delta x}\cdot \alpha_x + \frac{\Delta t}{\Delta y}\cdot \alpha_y \le 1$. %INSERIRE CENNO AL TEOREMA DI BARLES-SOUGANIDIS PER LA CONVERGENZA?
Moreover, as suggested in \cite{OS88}, this scheme can also be used to solve equation \eqref{eq:hamiltonian-2} by first splitting the Hamiltonian into two components 
$$
H(x,y,Dv,D^2 v)=H_0(x,y,Dv)+H_1(k(x,y),Dv),
$$
where $H_0$ includes the eikonal and advection term, whereas $H_1$ accounts for the dependence on the curvature $k(x,y)$, which is computed as
$$
k(x,y)=\nabla \cdot \left(\frac{Dv(x,y)}{|Dv(x,y)|}\right) = \frac{v_{xx} v_y^2 - 2v_xv_yv_{xy}+v_{yy}v_x^2}{{(v_x^2+v_y^2)}^{3/2}}.
$$
Then, it is enough to discretize the first-order part using the Lax-Friedrichs scheme and $H_1$ with simple centered finite difference approximations. In this case, the stability of the scheme requires a parabolic-type CFL condition of the form $\Delta t \sim \Delta x^2$.

\section{Pre-processing and post-processing}\label{pre_post_processing}
Before solving the evolution equations, we pre-processed the CT data by first applying a standard histogram equalization in order to highlight the relevant boundaries, and then a Gaussian filter with deviation $\sigma$ dependent on the noise of the given images. Further, before histogram equalization, we set all Hounsfield Unit (HU) values above 120 to 0 in order to reduce the relative weight of the bones. This process was executed slice-wise in a two dimensional manner. Finally, we created a mask for all HU values inside the interval $[5,120]$ to save the position of the non-muscular tissue, which is needed to improve the stability of the scheme for practical application.

The initialization of the scheme requires the user to select two points, one for each muscle. These points become the centers of the two spheres of radius $r=5\Delta x$ used as initial condition for the evolution, which is executed separately for each psoas. For all simulations the spatial discretization (the length of the voxel side) is set to $\Delta x = 0.1$, whereas the maximum number of iterations $N_{max}$ is heuristically chosen case by case, depending on the number of slices in the CT volume and the overall quality of the images.
% {\color{red}and the tolerance for the stopping criterion is set to $\tau = 1$; this helps the scheme avoid unwanted early stopping, particularly in the case of the second order model {\bf(eviterei di parlare del criterio di arresto, visto che non si usa mai, ed introdurrei il numero max di iterazioni)}}.

As previously remarked, clinical low-dose CT images as the ones acquired by PET/CT scanners, are usually noisy and present a lower level of resolution with respect to high-resolution CT acquisitions. Further, the complexity of the human body in the psoas region often implies the presence of sections where different organs or muscles, which have close HU values, have numerous contact points. Therefore, the definition of a clear boundary of the psoas becomes an intricate issue, which impacts the reliability of the numerical scheme in automatically segmenting the image. Specifically, in some cases it may happen that the wave-front leaves the correct boundary and starts detecting spurious objects outside the muscle of interest, thus making the usual stopping criterion ineffective, with unreliable final results. 
%{\color{red}This is why a maximum number of iterations $N_{max}$ is chosen case by case, depending on the number of slices in the CT and the overall quality of the images.}

In order to solve this issue, we devised a post-processing routine able to clean the result of (most of) the external objects and eventually give an acceptable final result. The procedure is composed of two steps and has a very low computational cost, which makes it efficient even if more than one iteration of the second step is needed. The routine works separately for each slice, but makes use of the information given by the closest slice in order to exploit the intrinsic continuity of the surrounding problem. We have summarized the procedure in \emph{Algorithm} 1, while more details on the implementation are given in the following items.
\begin{itemize}
\item {\bf\emph{First step.}}
The main aim of this step is to reduce the number of connected components (CCs) in each slice, detecting the one belonging to the psoas and discarding the remaining ones. To do so, the algorithm starts from the slice containing the center chosen by the user and select the correct CC by trivially excluding the ones that do not include the chosen point. Then, the remaining slices are scanned toward the top and the bottom, separately, always starting from the center. The main loop consists in first computing the centroid of the CC detected in the previous slice, and then, in the current slice, in selecting just the CC containing such point. The final result is a segmented muscle in which only the spurious objects directly connected to the boundary may be still present.
\item {\bf\emph{Second step.}} %This part is more intricate and has been proposed in two versions. 
%\emph{Version 1.} The first method has an easier implementation and more stable results, but it is not able to clean slices where spurious objects fall inside the region detected as correct in the previous slice. 
%This part of the procedure, which aims at cleaning all the remaining objects even if directly connected to the border, is more intricate and can be rather difficult to perform if all possible cases are considered. This is why we propose a simple implementation which has very stable results, but sometimes it is not able to clean slices where spurious objects fall inside the region detected as correct in the previous slice.
%More precisely, it is able to detect and delete the unwanted parts if the connected component is larger than the previous one, whereas it can have problems in detecting external objects if the region shrinks. In detail, the procedure assumes that the starting slice and the previous one, in the direction of the scanning, do not need the cleaning procedure. If this is not the case, a second iteration might be needed and the user has to select a new starting point.\\
%The initialization computes, for each point of the current slice, the minimum set distance w.r.t. the points belonging to the previous slice and stores it in a vector $tol$ of the same size as number of points in the current slice. 
In this second step the starting slice and the previous one, in the direction of the scanning, are assumed not to need the cleaning procedure. With this assumption, a vector $tol$ is initialized with length equal to $P$, the set of segmented voxels in the starting slice, where the $i-$th entry is the sum of minimum distance of the voxel $i \in P$ from the set of voxels $P_{prev}$ in the previous slice with a fixed constant (0.75).
Then the algorithm loops over the slice, setting $P_{prev} = P$ and first computing the difference set $P_{diff}=P-P_{prev}$, where $P$ and $P_{prev}$ are the sets of voxels in the current and previous slices, respectively; then, for each point $j \in P_{diff}$, the algorithm computes the set distance $dist(j,P_{prev})$ with respect to $P_{prev}$. This distance is obtained in correspondence of a voxel $i_{min} \in P_{prev}$. If $dist(j, P_{prev})>tol(i_{min})$, the voxel $j$ is excluded from $P$ and hence its removed from the segmentation. After all voxels in $P_{diff}$ have been scanned and removed if needed, a new tolerance vector $tol$ is computed using the current cleaned set of voxels $P$ and $P_{prev}$.
\end{itemize}

We point out that the second step, in its initialization phase, assumes that the starting slice and the previous one (according to the scanning direction) do not need any cleaning. If this is not the case, a second iteration of this step is needed, using a different starting point.
$$ $$
{\bf{Algorithm 1}}
\begin{description}
\item{{\bf{Input}}}: final solution $u$ of the segmentation method
\item{{\bf{step 1}}}: reduce the number of CCs in each slice
\begin{description}
\item{{\em{initialization:}}} start from the center and select the CC containing it 
\item{{\em{main loop:}}} scan the slice in both directions and keep the CC containing the centroid of the previous slice
\end{description}
\item{{\bf{step 2}}}: clean the remaining spurious objects: 
\begin{description}
\item{{\em{initialization:}}} Set a initial slice $S_i$. List all the voxels in the segmetation of that slice in the set $P$ and all the voxels in the segmetation of the previous slice $S_{i-1}$ in the set $P_{prev}$.
Initialize a vector $tol$ with length equal to cardinality of $P$. The value of the $i-$th entry of this vector is given by the sum minimum distance of voxel $i \in P$ from all the voxels in $P_{prev}$ and a fixed value (0.75). 
\item{{\em{main loop on the slices:}}} 
\begin{enumerate}
    \item set $P_{prev} = P$ 
    \item set $S_i =S_{i+1}$ and list all the voxels in the segmetation of slice $S_i$ in the set $P$
    \item compute the set difference  $\quad P_{diff} = P - P_{prev}$
    \item for each $j \in P_{diff}$: 
    \begin{itemize}
    \item compute its distance $dist(j,P_{prev}) $ w.r.t. the set $P_{prev}$ and save the point $i_{min}$ in $P_{prev}$ of minimum distance.
    \item remove the point $j$ from $P$ if $dist(j,P_{prev})>tol(i_{min})$
    \end{itemize}
     \item update the $tol$ vector, computing it as in the initialization step.
\end{enumerate}
\end{description}
\end{description}

\section{Numerical results}\label{numerical_results}
\begin{figure}
\centering
\includegraphics[width=0.8\textwidth]{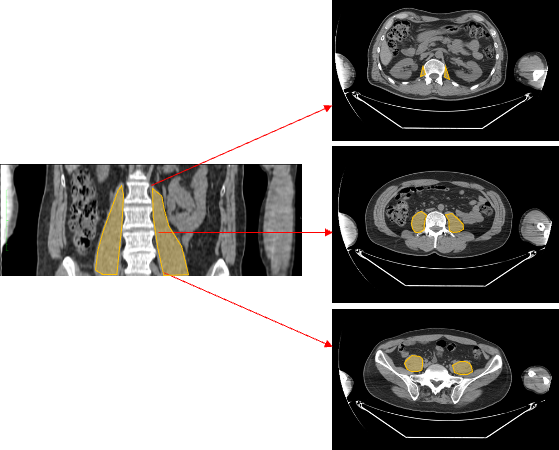}
\caption{Coronal view of the psoas muscles (in yellow) and the corresponding axial views at different levels of the abdomen.}\label{fig:Fig1_ok}
\end{figure}

The pipeline made of pre-processing, numerical segmentation, and post-processing has been applied to the CT images (see Figure \ref{fig:Fig1_ok}) acquired for nine subjects belonging to a retrospective data set of patients recruited at the IRCCS Ospedale Policlinico San Martino, Genova, Italy. All subjects were submitted to FDG-PET/CT and the study was performed according to the Declaration of Helsinki, Good Clinical Practice, and local ethics regulation. All enrolled patients signed a written informed consent at the time of FDG-PET/CT, econmpassing the use of anonymized data for retrospective research purposes.

For numerical segmentation we used the three numerical schemes introduced in Section 2, i.e., the finite difference scheme (\ref{eq:FD_2D})-(\ref{pippo-2}) applied to (\ref{eq_classic})  ('classical' from now on), and to the two geodesic schemes (\ref{eq:hamiltonian-2}) and (\ref{eq:hamiltonian-3}) ('GMFD 1ord' and 'GMFD 2ord' from now on, respectively). For both the two finite difference goedesic model schemes we set $\mu=1$ and $\nu=0.25$, while $\epsilon=0.05$ for the second order model. Further, 'GMFD 1ord' exploits the linear CFL condition $\Delta t = \lambda \Delta x$, with
\begin{equation}\label{lambda}
    \lambda = (3 \| g \| + \| D_x g \| + \| D_y g \| + \| D_z g\|)^{-1},
\end{equation}
and 'GMFD 2ord' utilizes the parabolic CFL condition $\Delta t = 0.25 \Delta x^2$.

%\subsubsection*{Test 1}

\begin{figure}
\centering
\begin{tabular}{ccc}
\includegraphics[width=0.25\textwidth]{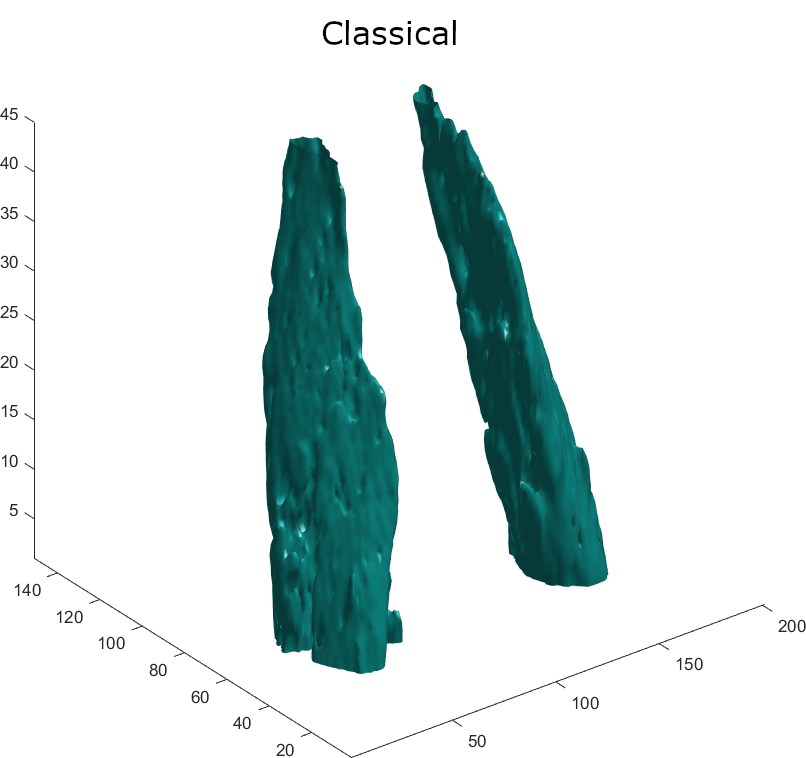} &
\includegraphics[width=0.25\textwidth]{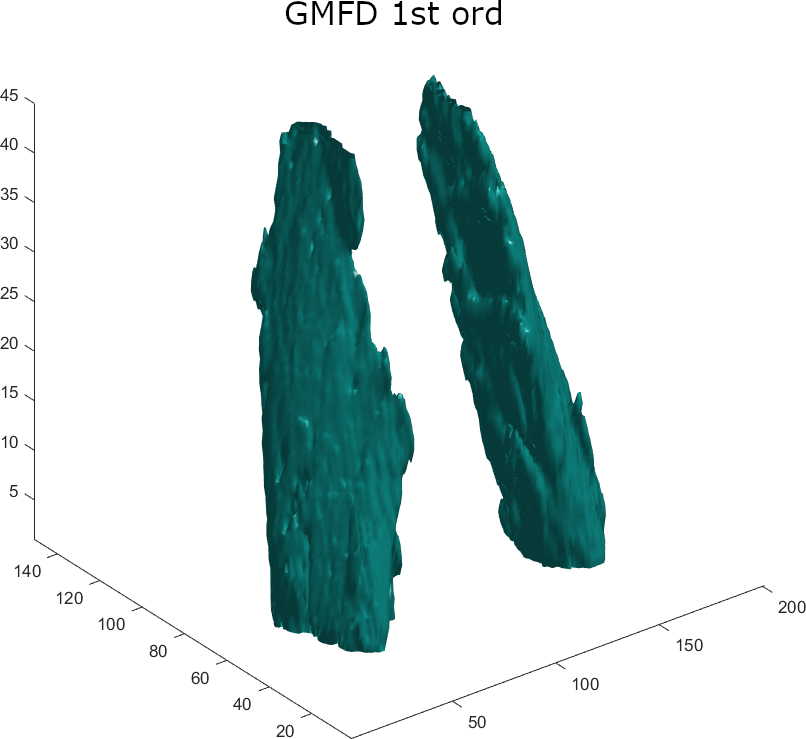} &
\includegraphics[width=0.25\textwidth]{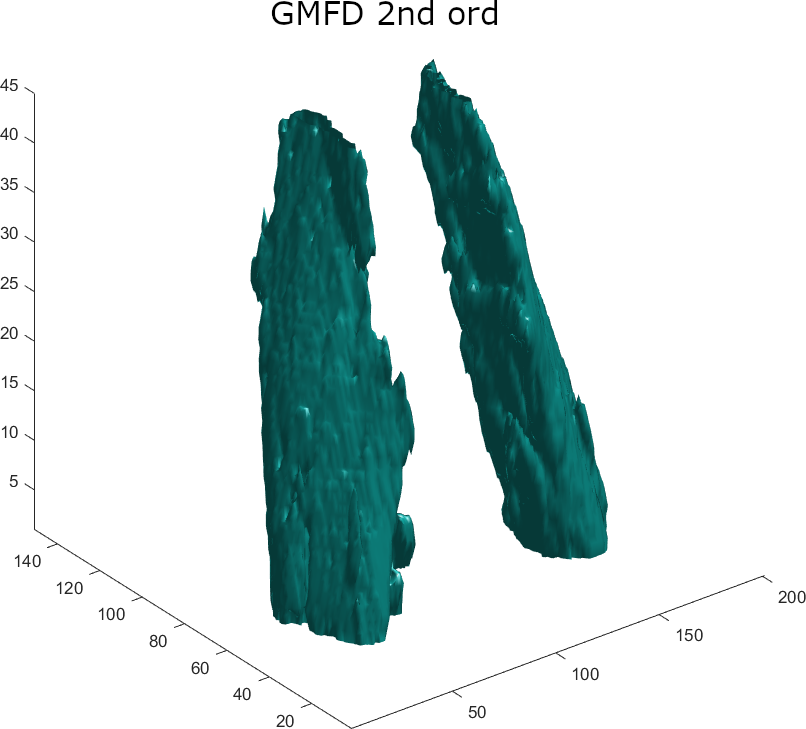} 
\\
\includegraphics[width=0.25\textwidth]{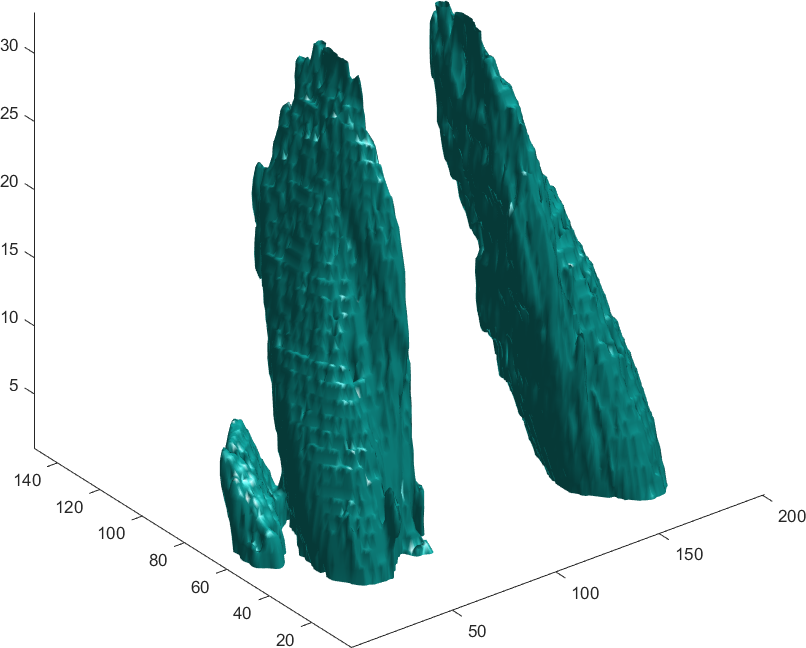} &
\includegraphics[width=0.25\textwidth]{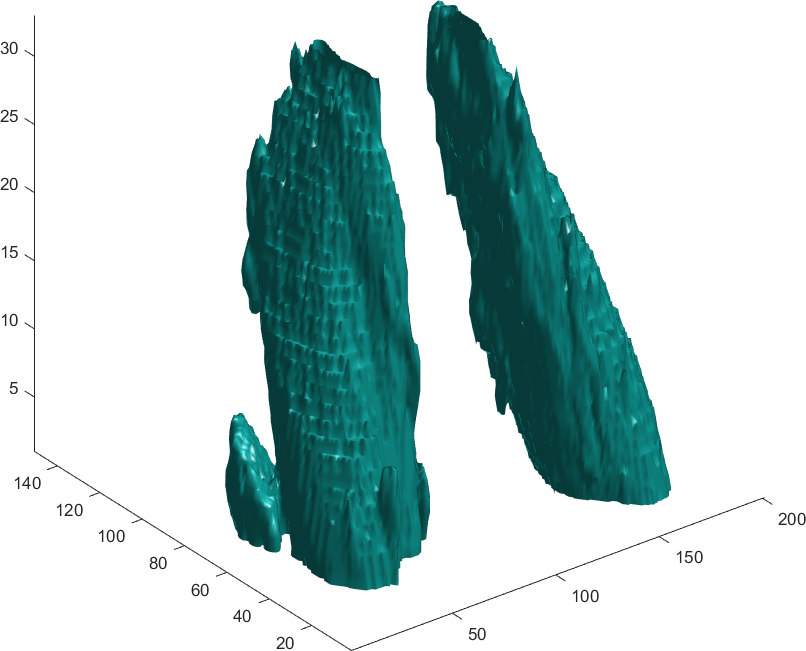} &
\includegraphics[width=0.25\textwidth]{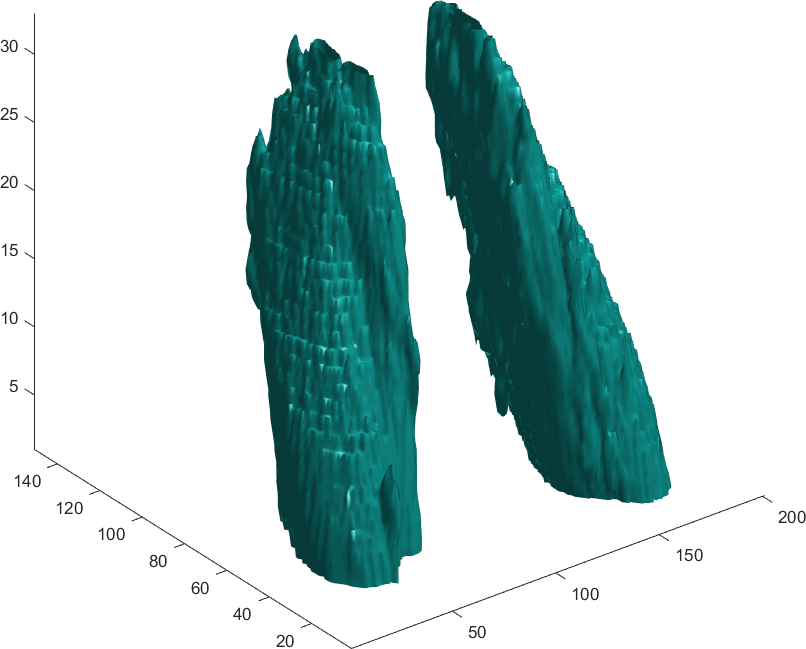} \\
\includegraphics[width=0.25\textwidth]{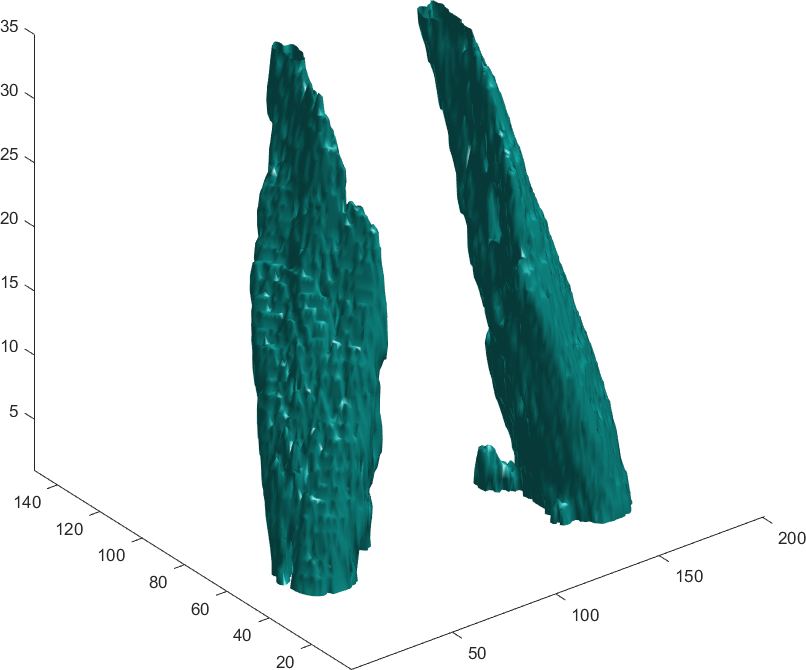} &
\includegraphics[width=0.25\textwidth]{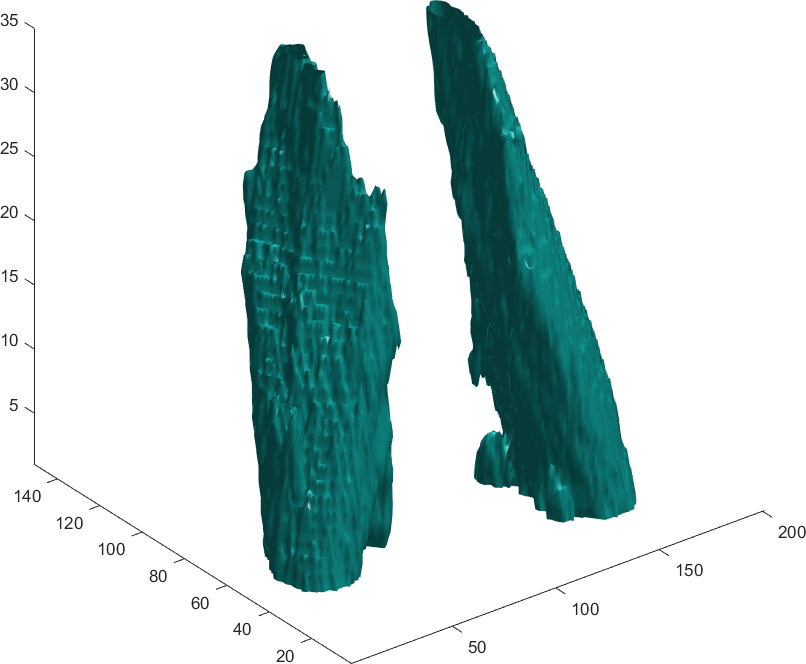} &
\includegraphics[width=0.25\textwidth]{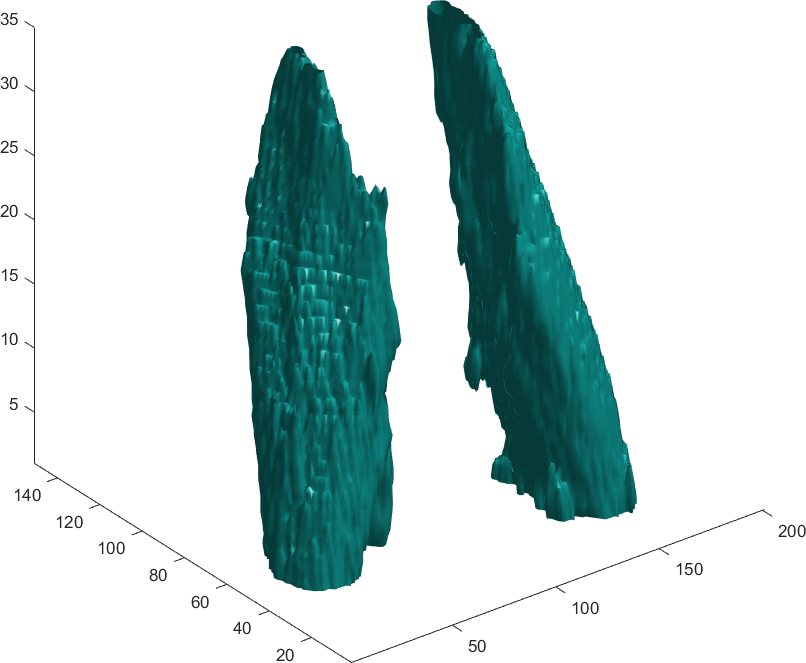} \\
\includegraphics[width=0.25\textwidth]{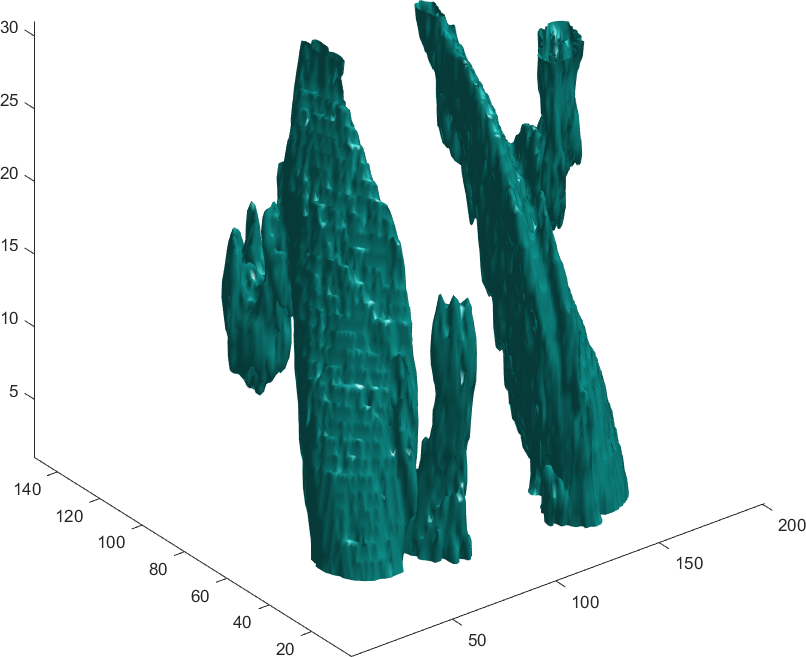} &
\includegraphics[width=0.25\textwidth]{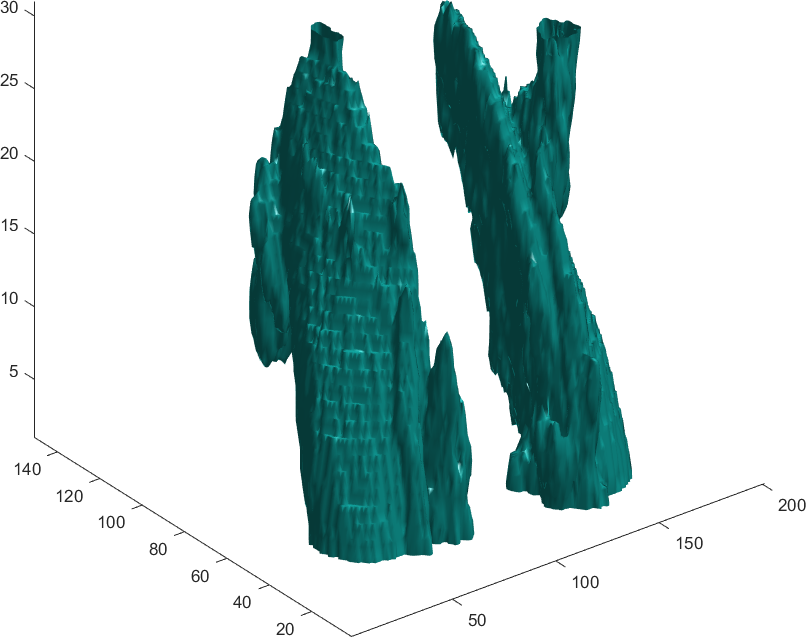} &
\includegraphics[width=0.25\textwidth]{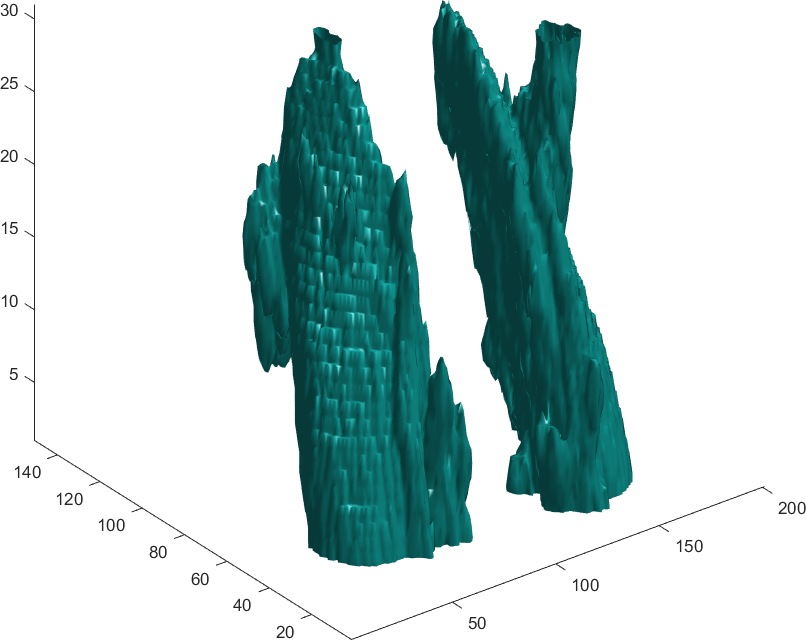}
\\
\includegraphics[width=0.25\textwidth]{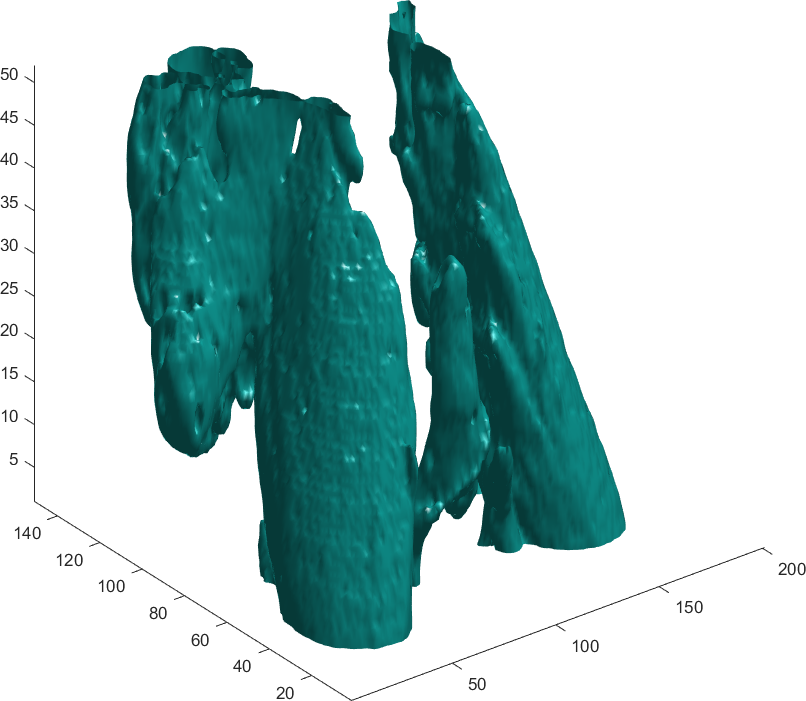} &
\includegraphics[width=0.25\textwidth]{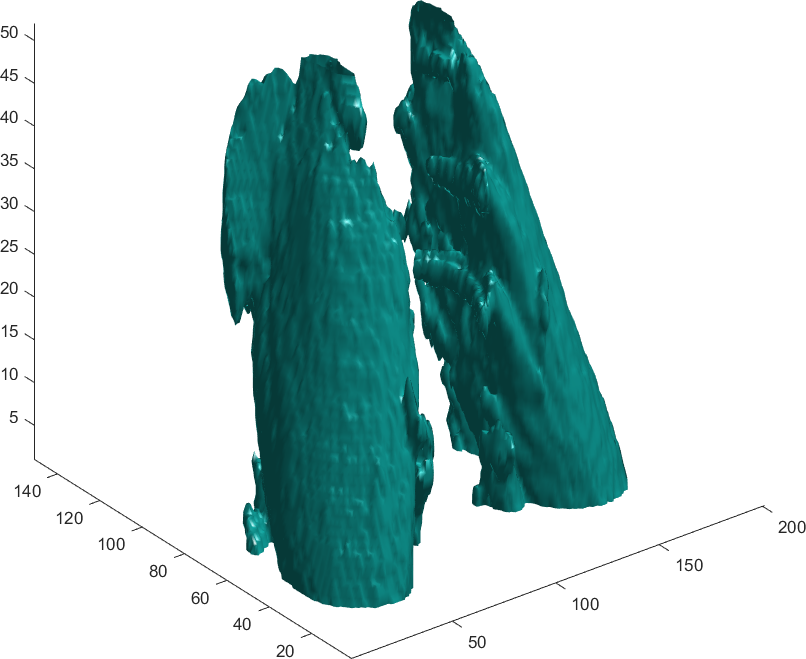} &
\includegraphics[width=0.25\textwidth]{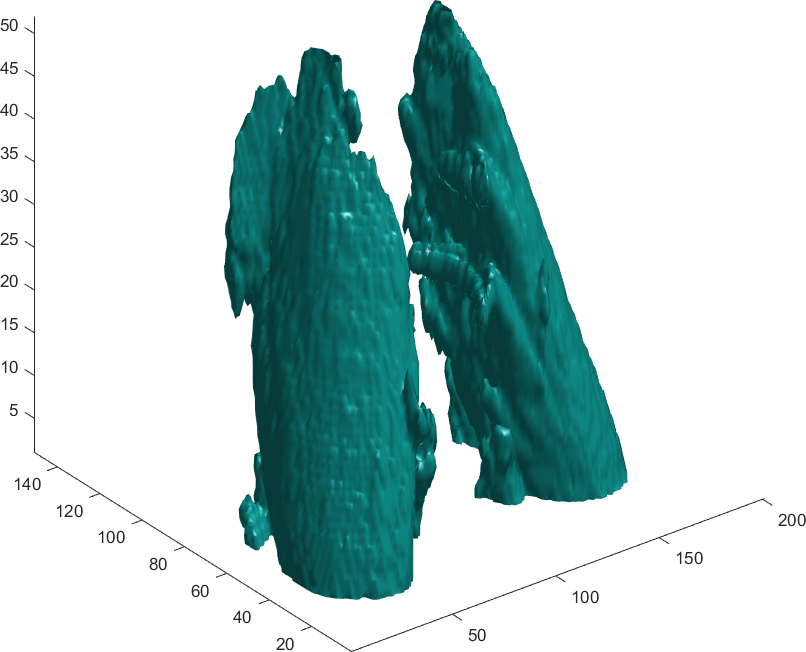}
\\
\includegraphics[width=0.25\textwidth]{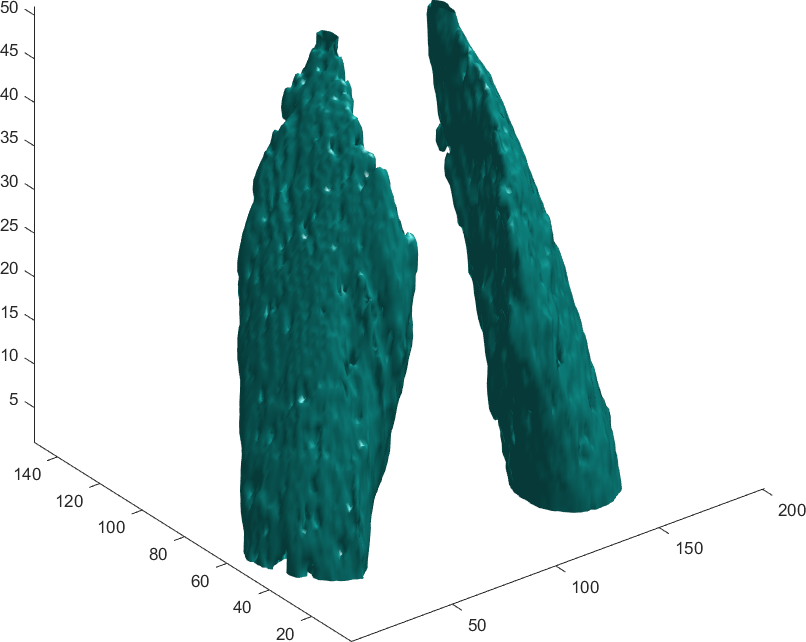} &
\includegraphics[width=0.25\textwidth]{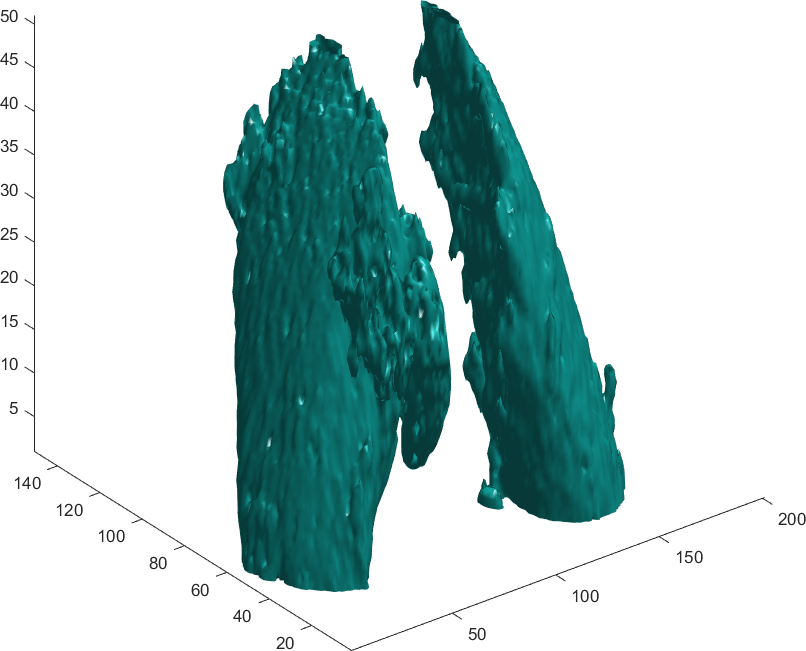} &
\includegraphics[width=0.25\textwidth]{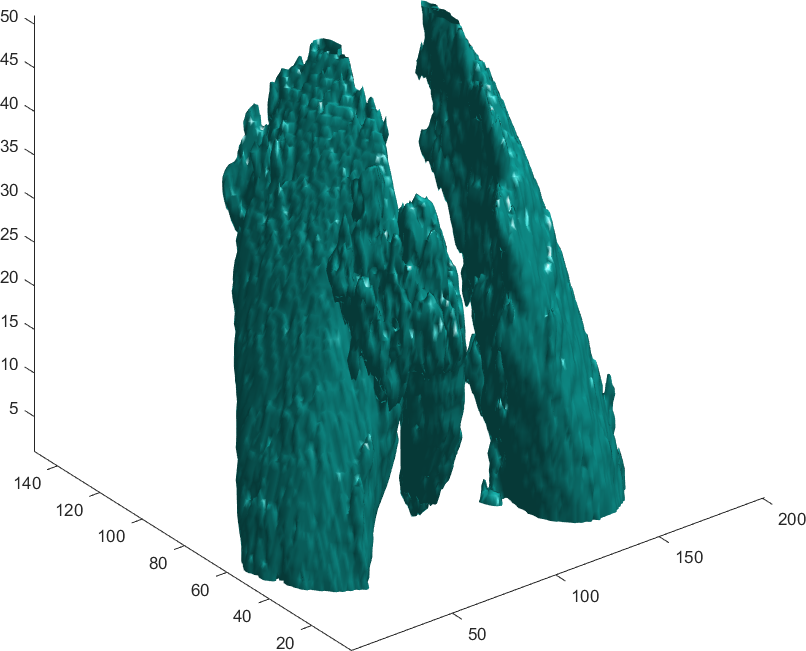}
\end{tabular}
\caption{Segmentations provided by the classical evolution model (\ref{eq_classic}) ('classical', first column), the first-order geodesic model (\ref{eq:hamiltonian-3}) ('GMFD 1ord', second column), and the second order geodesic model (\ref{eq:hamiltonian-2}) ('GMFD 2ord', third column). First row: subject CTTO02 (45 slices) with $p=2$ and $\sigma=0$. Second row: subject PZGE01 (33 slices) with $p=2$ and $\sigma=2$. Third row: subject PZGE02 (35 slices) with $p=2$ and $\sigma=2$. Fourth row: subject PZGE03 (31 slices) with $p=2$ and $\sigma=2$. Fifth row: subject  PZTO03 (52 slices) with $p=2$ and $\sigma=0$. Sixth row: subject PZTO04 (51 slices) with $p=2$ and $\sigma=1$.}\label{fig:results-1}
\end{figure}

\begin{figure}
\centering
\includegraphics[width=0.30\textwidth]{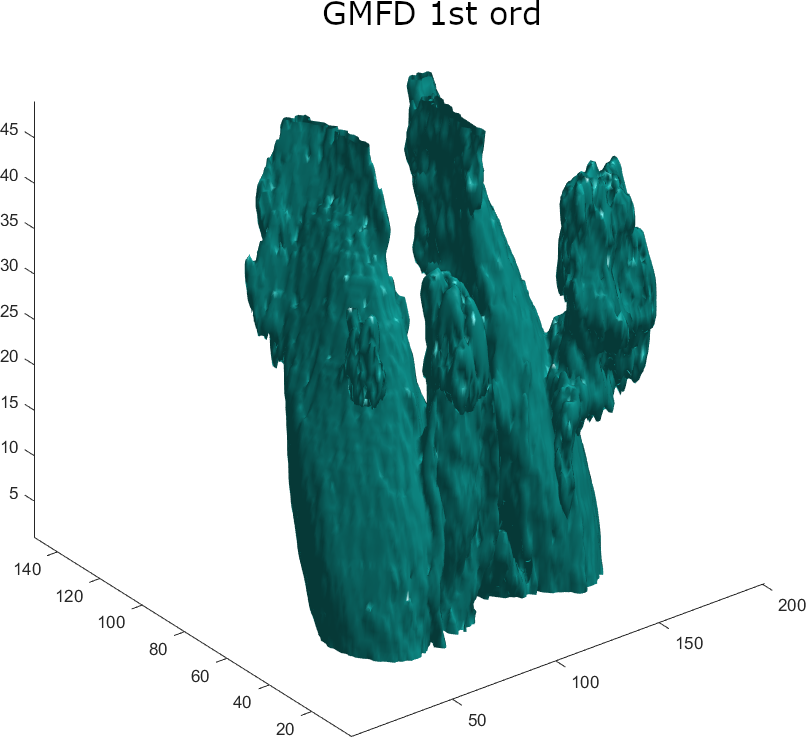}
\includegraphics[width=0.30\textwidth]{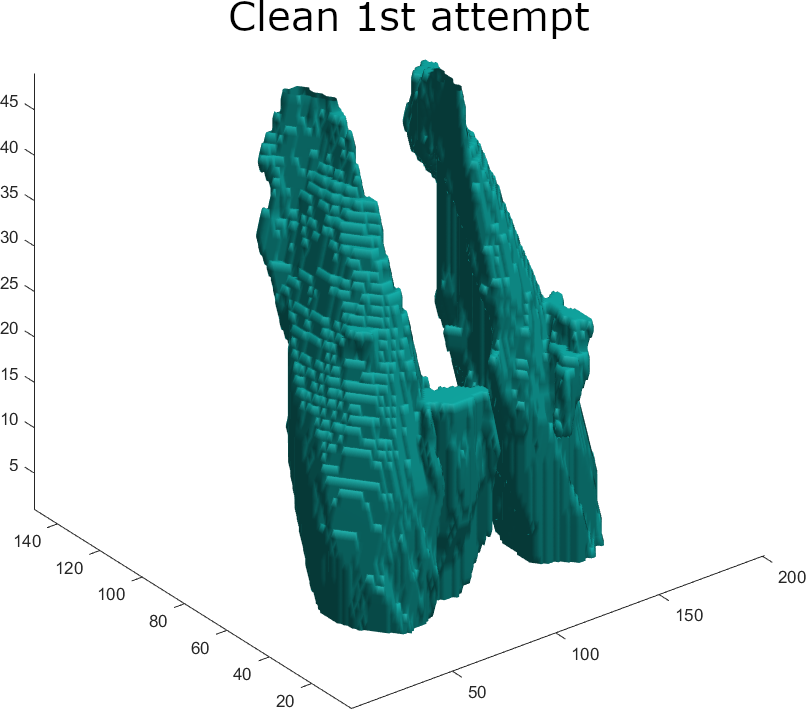}
\includegraphics[width=0.30\textwidth]{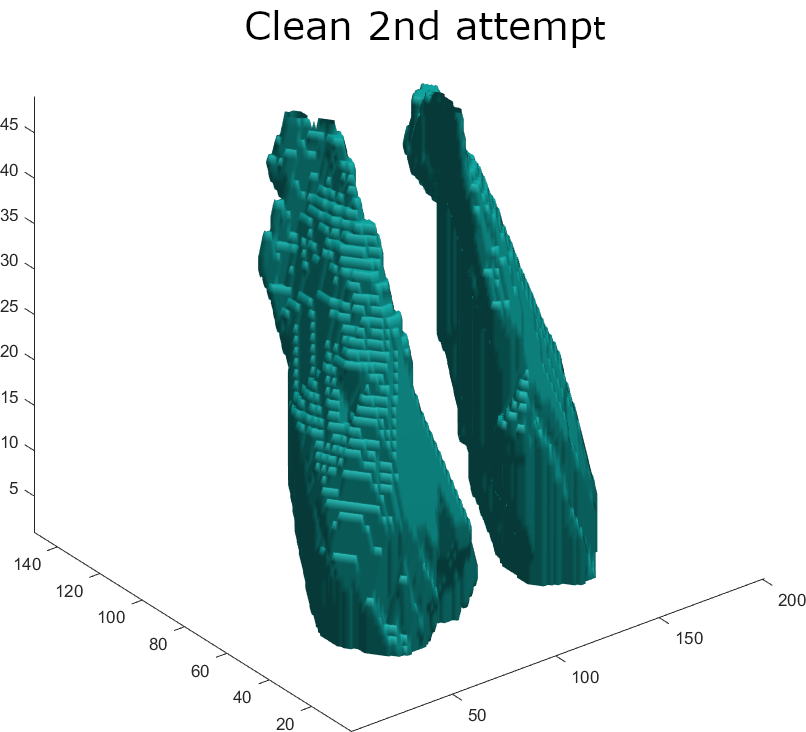} \\
\includegraphics[width=0.30\textwidth]{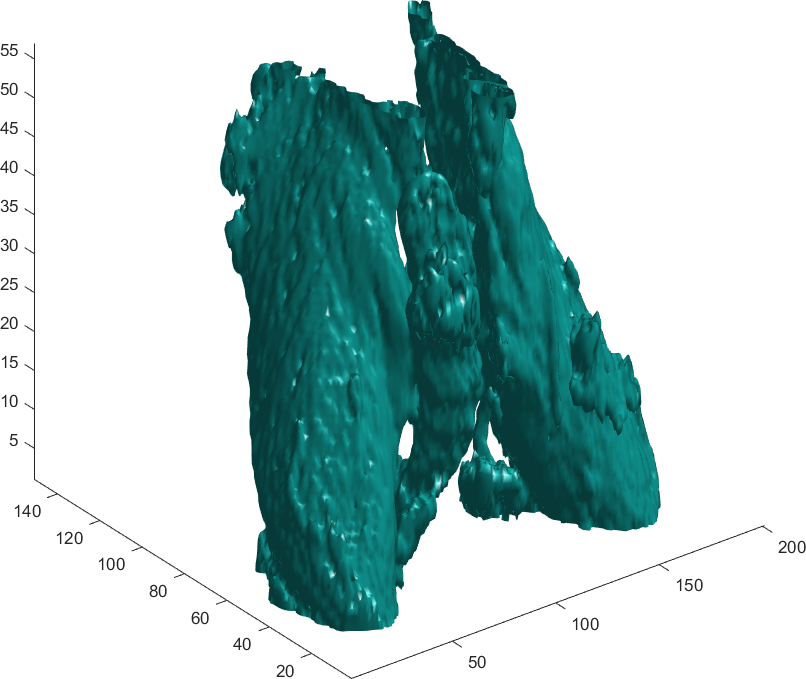}
\includegraphics[width=0.30\textwidth]{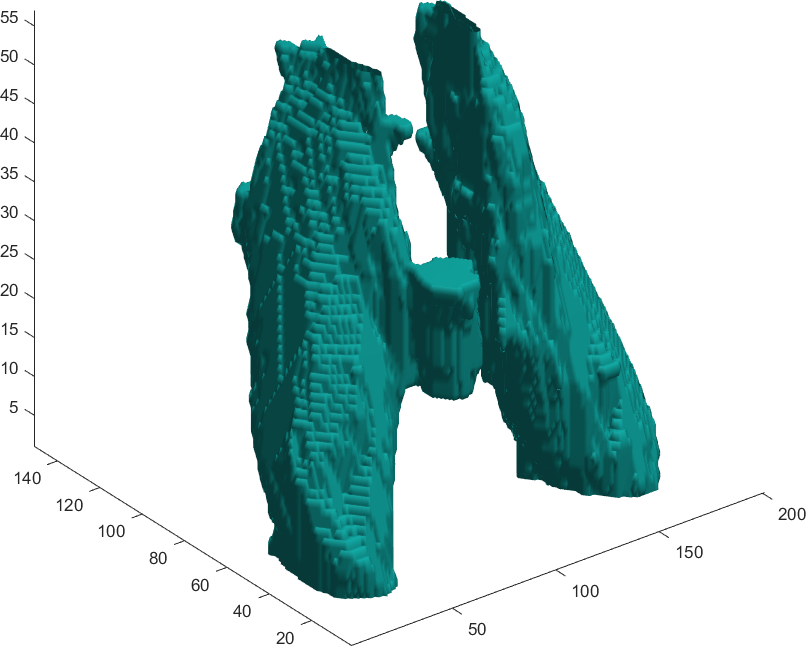}
\includegraphics[width=0.30\textwidth]{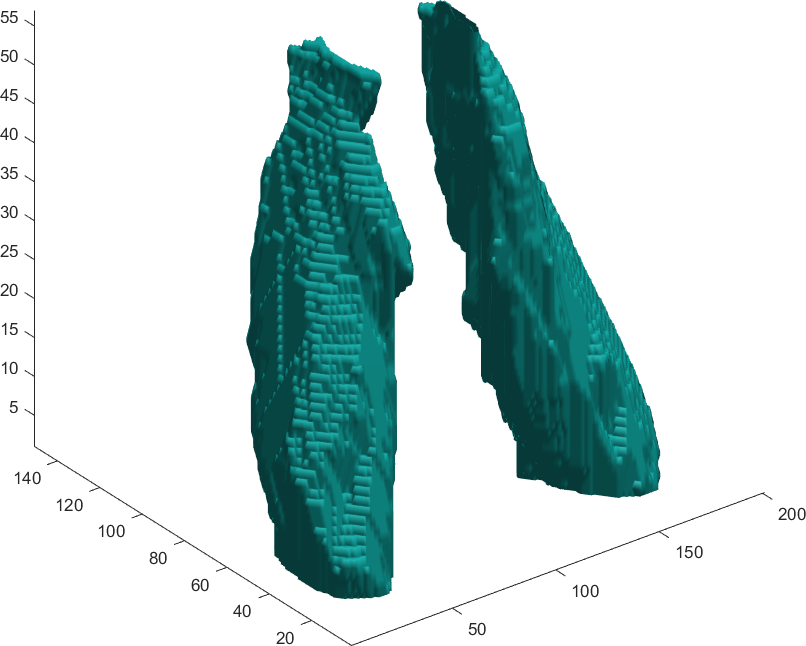} \\
\includegraphics[width=0.30\textwidth]{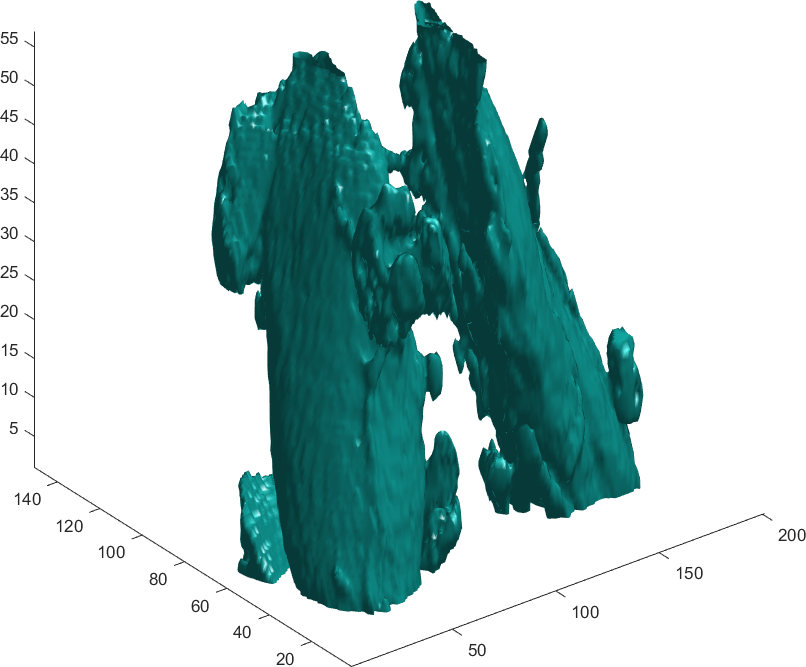}
\includegraphics[width=0.30\textwidth]{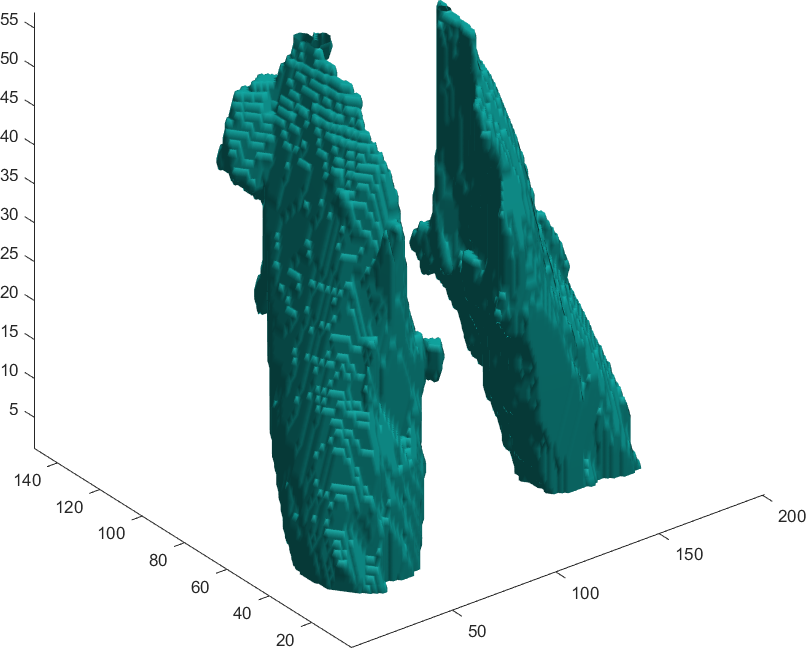}
\includegraphics[width=0.30\textwidth]{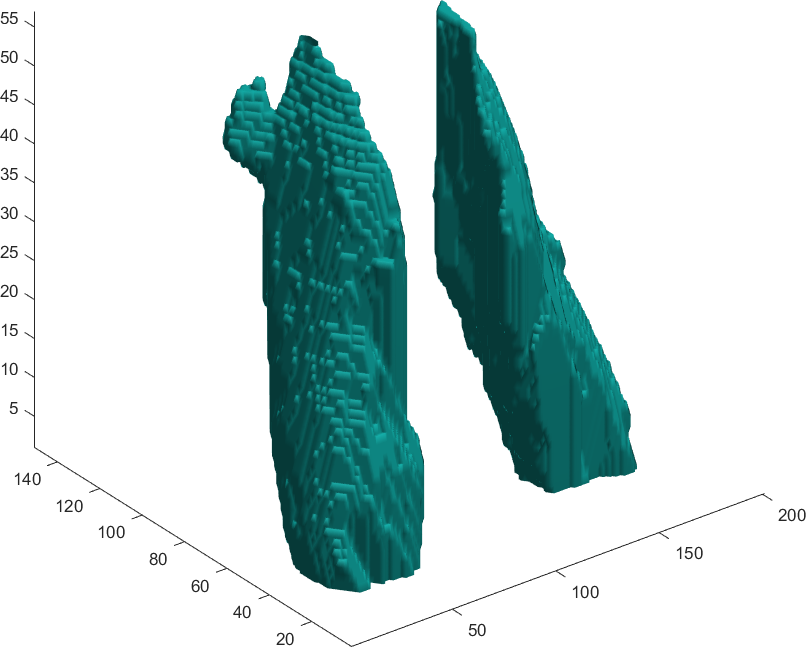}
\caption{\small{Post-processing pipeline for the patients CTTO01, CTTO03, and PZTO01 when applied to the outcomes of 'GMFD 1ord'. First column: the starting data. Second colum: the outcome of the first step (starting slice: $20$ for CTTO01, $44$ for CTTO03, $25$ for PZTO01). Third column: the outcome of the second step.}\label{fig:cl_test1_2}}
\end{figure}

Table \ref{table:results-1} contains the number of iterations and CPU time employed to realize the segmentation by the three numerical schemes when applied to the volumes of all nine patients. Figure \ref{fig:results-1} contains the results of the segmentation provided by the three methods in the case of six subjects and Figure \ref{fig:cl_test1_2} describes the performances of the post-processing algorithm in the case of the remaining three subjects. Specifically, the three cases in this latter figure are the ones for which the outcomes of the numerical schemes present a particularly significant amount of artifacts. However, the post-processing step can be applied to all outcomes, including the ones described in Figure \ref{fig:results-1}.

\begin{table}[h!]
\begin{tabular}{c | c | c | c | c | c | c}
method & patient & $N_{max}$ & CPU & $p$ & $\sigma$ & slices \\
\hline\hline
classical & CTTO1 & 400 & 178.48 & 2 & 1 & 49 \\
GMFD 1ord & CTTO1 & 400 & 197.59 & 2 & 1 & 49 \\
GMFD 2ord & CTTO1 & 1900 & 1955.84 & 2 & 1 & 49 \\
\hline
classical & CTTO02 & 350 & 133.19 & 2 & 0 & 45 \\
GMFD 1ord & CTTO02 & 350 & 147.56 & 2 & 0 & 45 \\
GMFD 2ord & CTTO02 & 2000 & 1870.16 & 2 & 0 & 45 \\
\hline
classical & CTTO03 & 450 & 232.91 & 2 & 1& 57 \\
GMFD 1ord & CTTO03 & 450 & 250.65 & 2 & 1 & 57 \\
GMFD 2ord & CTTO03 & 2700 & 3298.88 & 2 & 1 & 57 \\
\hline
classical & PZGE01 & 700 & 229.96 & 2 & 2 & 33 \\
GMFD 1ord & PZGE01 & 700 & 234.94 & 2 & 2 & 33 \\
GMFD 2ord & PZGE01 & 3500 & 2654.22 & 2 & 2 & 33 \\
\hline
classical & PZGE02 & 700 & 232.67 & 2 & 2 & 35 \\
GMFD 1ord & PZGE02 & 700 & 246.20 & 2 & 2 & 35 \\
GMFD 2ord & PZGE02 & 3500 & 2787.62 & 2 & 2 & 35 \\
\hline
classical & PZGE03 & 700 & 207.53 & 2 & 2 & 31 \\
GMFD 1ord & PZGE03 & 700 & 224.22 & 2 & 2 & 31 \\
GMFD 2ord & PZGE03 & 3500 & 2499.52 & 2 & 2 & 31 \\
\hline
classical & PZTO01 & 500 & 258.00 & 2 & 0& 57 \\
GMFD 1ord & PZTO01 & 500 & 275.68 & 2 & 0 & 57 \\
GMFD 2ord & PZTO01 & 2500 & 2991.82 & 2 & 0  & 57 \\
\hline
classical & PZTO03 & 500 & 229.58 & 2 & 0 & 52 \\
GMFD 1ord & PZTO03 & 500 & 250.62 & 2 & 0 & 52 \\
GMFD 2ord & PZTO03 & 2500 & 2758.94 & 2 & 0 & 52 \\
\hline
classical & PZTO04 & 450 & 207.79 & 2 & 1 & 51 \\
GMFD 1ord & PZTO04 & 450 & 230.97 & 2 & 1 & 51 \\
GMFD 2ord & PZTO04 & 2300 & 2540.63 & 2 & 1 & 51 \\
\hline\hline
\end{tabular}
\caption{Number of iterations, CPU time, values of $p$ and $\sigma$, and number of slices for the three segmentation algorithms for the nine CT volumes considered in this study.}\label{table:results-1}
\end{table}

We compared the segmentation results we obtained with the three methods described in the previous sections with the ones provided by manually drawn profiles within the OsiriX
software package \citep{OsiriX} and by the 3D Slicer extension TotalSegmentator (TS) \citep{TotalSegmentator} used in both fast and full-resolution modes.
In order to evaluate the agreement between the manually drawn psoas profiles and the segmentations obtained with the three numerical schemes and the TotalSegmentator tool we used four metrics: the Dice similarity coefficient \citep{Dice}, the Jaccard index \citep{Jaccard}, the Hausdorff distance \citep{hausdorff1914grundzuge} and the Average Symmetric Surface Distance (ASSD) \citep{ASSD}. The Dice similarity coefficient and the Jaccard index vary between $0$ and $1$ where $0$ represents an empty intersection between the two segmentation while $1$ a perfect overlap. The Hausdorff distance and ASSD are indeed distances, so smaller values correspond to closer objects. Specifically, the Hausdorff distance is the greatest distance from a point in one segmentation mask to the closest point in the other mask, while ASSD measures the average of all distances from points on the boundary of the first mask to the boundary of the second one, and viceversa. 
In Table \ref{tab:distanze} we show the average values of these metrics computed across the nine subjects, together with the corresponding standard deviations.

\begin{table}[h!]
\begin{tabular}{c | c | c | c | c | c }
metric & 'classical' & 'GMFD 1ord & 'GMFD 2ord' & TS-FM & TS-FRM\\
\hline\hline
Dice	&	0.86	$\pm$	0.03	&	0.87	$\pm$	0.02	&	0.87	$\pm$	0.03	&	0.83	$\pm$	0.02	&	0.85	$\pm$	0.05	\\
Jaccard	&	0.76	$\pm$	0.05	&	0.77	$\pm$	0.04	&	0.77	$\pm$	0.04	&	0.72	$\pm$	0.03	&	0.75	$\pm$	0.07	\\
Hausdorff	&	12.15	$\pm$	3.85	&	12.95	$\pm$	2.96	&	12.64	$\pm$	3.76	&	29.18	$\pm$	25.95	&	29.55	$\pm$	15.58	\\
ASSD	&	0.40	$\pm$	0.10	&	0.40	$\pm$	0.04	&	0.42	$\pm$	0.08	&	0.64	$\pm$	0.35	&	0.66	$\pm$	0.38	\\
\hline\hline
\end{tabular}
\caption{Reliability assessment for the three three-dimensional numerical schemes when compared to manual segmentation. The metrics utilized for this analysis are the Dice coefficient (optimal value: $1$), the Jaccard coefficient (optimal value: $1$), and the Hausdorff and ASSD distances, respectively. The last two columns report the metrics values in the case of two standard software tools utilized in medical imaging, i.e., the TotalSegmentator 3D extension of Slicer in both fast and full-resolution mode (TS-FM and TS-FRM, respectively).}\label{tab:distanze}
\end{table}

\section{Comments and conclusions}\label{conclusions}
This paper has proposed a 3D approach to the segmentation of the psoas muscle based on level set models. These numerical schemes provide a higher degree of automation with respect to the approach utilized in \cite{bauckneht2020spinal}, which is notably heuristic, and a higher degree of generality with respect to the approach utilized in \cite{bauckneht2022opportunistic}, whose performances rely on the availability of prior geometrical models for the psoas.

For all numerical tests the geodesic models provide stable results, whereas the classical scheme may produce over-resolved shape, as it is particularly evident in the cases of patient CTTO02 (Figure \ref{fig:results-1}, first column, first row) and patient PZGE02 (Figure \ref{fig:results-1}, first column, third row). This effect is not evident for the two geodesic models, probably because the transport term in \eqref{eq:hamiltonian-2}-\eqref{eq:hamiltonian-3} pushes the front toward the external boundary of the muscle. Further, this same transport term usually helps the front stick to the right boundary, reducing the overflow outside the muscle, which happens when some complex structure is present in some portion of the data, as shown in Figure \ref{fig:results-1}, fifth row. As far as the reliability of the segmentation is concerned, the average Hausdorff and ASSD distances computed between the masks produced by our schemes and the manually identified ground truth are, respectively, more than two times and more than 1.5 times lower than the distances obtained using the TotalSegmentator masks. This fact highlights that, even if the Jaccard and the Dice indicate comparable overlap, our segmentation masks present less outlier. 

The first order version of the geodesic model has performances with comparable reliability to the ones characterizing the second order model. Indeed, the metrics and distances used in Table \ref{tab:distanze} have similar mean values and uncertainties. These same uncertainties are clearly smaller than the ones associated to the classical model. Further, the computational cost required by 'GMFD 1ord' to achieve its metrics and distances values is one order of magnitude smaller than the one typical of the second order model (see Table \ref{table:results-1}). The higher computational burden of 'GMFD 2ord' is mainly due to the parabolic CFL condition, which slows down the evolution and forces the algorithm to utilize a greater number of iterations to reach convergence
.
It is clear that the simple post-processing procedure devised to clean the results provided by the numerical schemes is able to remove all spurious objects in most cases, still keeping the volume of the segmented muscle mostly untouched with respect to the results of segmentation. From a heuristic viewpoint, we could test that, if the outcome given by the scheme is already accurate on its own, one iteration of the procedure is enough to remove most artifacts. On the other hand, if the spurious objects are many and, in particular, some of them intersect the slice containing the center of the initial condition, a second run with a (possibly) different starting point is necessary. In any case, if the latter point is chosen correctly, the final result can be satisfactory and only in few cases some small objects still persist. 

We finally observe that higher order schemes, as defined for example in \cite{FPT20}, could be considered and adapted to the modified model \eqref{eq:hamiltonian-2} to increase the accuracy of the segmentation at the expense of longer computational times, but this goes beyond the scope of the present work and will be considered in future research.

\section*{Acknowledgements}
CC and MP acknowledge the financial support of the "Hub Life Science – Digital Health (LSH-DH) PNC-E3-2022-23683267 - Progetto DHEAL-COM – CUP: D33C22001980001", granted by the Italian Ministero della Salute within the framework of the Piano Nazionale Complementare to the "PNRR Ecosistema Innovativo della Salute - Codice univoco investimento: PNC-E.3"

GP and MP dedicate this paper to the memory of Maurizio Falcone.

\bibliographystyle{aa}
\bibliography{sn-bibliography}

\end{document}